\documentstyle[12pt,amssymb]{article}

\newtheorem{Theo}{Theorem}
\newtheorem{ex}{Example}

\mathchardef\za="710B  
\mathchardef\zb="710C  
\mathchardef\zg="710D  
\mathchardef\zd="710E  
\mathchardef\zve="710F 
\mathchardef\zz="7110  
\mathchardef\zh="7111  
\mathchardef\zvy="7112 
\mathchardef\zi="7113  
\mathchardef\zk="7114  
\mathchardef\zl="7115  
\mathchardef\zm="7116  
\mathchardef\zn="7117  
\mathchardef\zx="7118  
\mathchardef\zp="7119  
\mathchardef\zr="711A  
\mathchardef\zs="711B  
\mathchardef\zt="711C  
\mathchardef\zu="711D  
\mathchardef\zvf="711E 
\mathchardef\zq="711F  
\mathchardef\zc="7120  
\mathchardef\zw="7121  
\mathchardef\ze="7122  
\mathchardef\zy="7123  
\mathchardef\zf="7124  
\mathchardef\zvr="7125 
\mathchardef\zvs="7126 
\mathchardef\zf="7127  
\mathchardef\zG="7000  
\mathchardef\zD="7001  
\mathchardef\zY="7002  
\mathchardef\zL="7003  
\mathchardef\zX="7004  
\mathchardef\zP="7005  
\mathchardef\zS="7006  
\mathchardef\zU="7007  
\mathchardef\zF="7008  
\mathchardef\zW="700A  

\newcommand{\be}{\begin{equation}}
\newcommand{\ee}{\end{equation}}
\newcommand{\lra}{\longrightarrow}
\newcommand{\bea}{\begin{eqnarray}}
\newcommand{\eea}{\end{eqnarray}}
\newcommand{\beas}{\begin{eqnarray*}}
\newcommand{\eeas}{\end{eqnarray*}}

\newcommand{\R}{{\Bbb R}}
\newcommand{\C}{{\Bbb C}}
\newcommand{\SL}{SL(2,\C)}

\newcommand{\SU}{SU(2)}
\newcommand{\su}{su(2)}
\newcommand{\SB}{SB(2,\C)}
\newcommand{\Sb}{sb(2,\C)}
\newcommand{\G}{{\cal G}}
\newcommand{\D}{{\rm d}}

\newcommand{\we}{\wedge}
\newcommand{\nn}{\nonumber}
\newcommand{\ot}{\otimes}
\newcommand{\s}{{\textstyle *}}
\newcommand{\ts}{T^\s}
\newcommand{\pa}{\partial}
\newcommand{\ti}{\times}
\begin{document}
\begin{center}
{\Large\bf Completely integrable systems:\\ a generalization}
\vskip 1cm
{D. V. Alekseevsky\footnote
{Gen. Antonova 2 kv 99, 117279 Moscow B-279, Russia},
J. Grabowski\footnote{Institute of Mathematics, Warsaw University,
ul. Banacha 2, 02-097 Warszawa, Poland; and \\
Mathematical Institute, Polish Academy of Sciences, ul.
\'Sniadeckich 8, P.O. Box 137, 00-950 Warszawa, Poland\\
{\it e-mail:} jagrab@mimuw.edu.pl},
G. Marmo\footnote{Dipartimento di Scienze Fisiche, Universit\`a di Napoli,
Mostra d'Oltremare, Pad. 19, 80125 Napoli, Italy\\
{\it e-mail:} gimarmo@na.infn.it}, and
P. W. Michor\footnote{Institut f\"ur Mathematik, Universit\"at Wien,
Strudlhofgasse 4, 1090 Wien, Austria; and \\Erwin Schr\"odinger International
Institute of Mathematical Physics, Bolzmanngasse 9, 1090 Wien, Austria\\
{\it e-mail:} michor@esi.ac.at}}
\vskip .5cm
{\it Erwin Schr\"odinger International Institute \\
of Mathematical Physics, Wien,Austria}
\end{center}
\date{\ }
\centerline{\bf Abstract}
We present a slight generalization of the notion of completely integrable
systems to get them being integrable by quadratures. We use this
generalization
to integrate dynamical systems on double Lie groups.
\bigskip\noindent
\setcounter{equation}{0}
\section{Introduction}
A Hamiltonian system on a $2n$-dimensional symplectic manifold $M$
is said to be {\it completely integrable} if it has $n$ first
integrals in involution,
which are functionally independent on some open dense submanifold
of $M$. This definition of a completely integrable system is
usually found, with some minor variants, in any modern text on
symplectic mechanics \cite{AM, Ar, LM, MSS, Th}.
\par
Starting with this definition, one uses the so called
Liouville-Arnold theorem to introduce action-angle variables and
write the Hamiltonian system in the form
\bea
\label{x} \dot I^k&=&0,\\
\label{y} \dot \zvf_k&=&\frac{\pa H}{\pa I_k}=\zn_k(I),
\eea
where $k\in\{ 1,\dots ,n\}.$
The corresponding flow is given by
\bea\label{1}
I^k(t)&=&I^k(0),\\
\label{2} \zvf_k(t)&=&\zvf_k(0)+\zn_kt.
\eea
The main interest in completely integrable systems relies on the fact that
they can be integrated by quadratures.
\par
It is clear, however, that even if $\ \zn_k\D I^k\ $ is not an
exact (or even a closed) 1-form, as long as $\dot\zn_k=0,$ the
system can always be integrated by quadratures.
\par
In this letter, we would like to take up this remark from the view
point of Lie groups and their cotangent bundles, as well as double
Lie groups, plying the
role of deformed cotangent bundles, to show how the notion of a completely
integrable system can be generalized by relaxing the property of being
Hamiltonian and retaining only that it has enough constants of the motion
(first integrals) to warrant it to be integrable by quadratures.
More precisely, we define some class of dynamical systems on a
Poisson manifold $(M,\zL)$, defined by a vector field $\zG
=\zL(\zh)$, where $\zh$ is a 1-form, which are integrable by
quadratures. Here $M$ is the cotangent bundle of a Lie group or
its appropriate deformation.

\setcounter{equation}{0}
\section{A universal model for completely integrable system}
If we consider the abelian Lie group $\R^n,$ we can construct a
Hamiltonian action of $\R^n$ on $\ts\R^n$ induced by the group
addition:
\be
\R^n\ti\ts\R^n\lra\ts\R^n.
\ee
This can be generalized to the Hamiltonian action
\be
\R^n\ti\ts(\R^k\ti T^{n-k})\lra\ts(\R^k\ti T^{n-k}),
\ee
of $\R^n$, where $T^m$ stands for the $m$-dimensional torus, and
reduces to
$\R^n\ti\ts T^n $ or $T^n\ti\ts T^n,$ when $k=0.$
\par
By using the standard symplectic structure on $\ts\R^n$, we find
the  momentum map $\zm:\ts\R^n\lra (\R^n)^\s,\ (q,p)\mapsto p$,
induced by the natural action of $\R^n$ on itself via
translations, which is a
Poisson map if $(\R^n)^\s$ is endowed with the (trivial) natural
Poisson structure of the dual of a Lie algebra.
It is now clear that any function on $(\R^n)^\s$, when pulled back to
$\ts\R^n$
or $\ts T^n$, gives rise to a Hamiltonian system which is completely
integrable
(in the Liouville sense).
Because the level sets of this function carry on the action of $\R^n$, the
completely integrable system gives rise to a one-dimensional subgroup of the
action of $\R^n$ on the given level set. The specific subgroup will, however,
depend on the particular level set, i.e. the `frequencies' are first
integrals.
The property of being integrable by quadratures is captured by the fact that
it
is a subgroup of the $\R^n$-action on each level set.
\par
It is now clear, how we can preserve this property, while giving up the
requirement that our system is Hamiltonian. We can indeed consider any 1-form
$\zh$ on $(\R^n)^\s$ and pull it back to $\ts\R^n$ or $\ts T^n$, then
associated vector field $\zG_\zh=\zL_0(\zm^\s(\zh))$, where $\zL_0$
is the canonical Poisson structure in the cotangent bundle, is no more
Hamiltonian, but it is still integrable by
quadratures. In action-angle variables, if $\zh=\zn_k\D I^k$ is the 1-form on
$(\R^n)^\s$, the associated equations of motion on $\ts T^n$ will be
\bea
\dot I^k&=&0,\\
\dot \zvf_k&=&\zn_k,
\eea
with $\dot\zn_k=0,$ therefore the flow will be as in (\ref{1}), (\ref{2}),
even though
\be
\frac{\pa\zn_k}{\pa I^j}\ne\frac{\pa\zn_j}{\pa I^k}.
\ee
We can now generalize this construction to any Lie group $G$. We consider the
Hamiltonian action
\be
G\ti\ts G\lra\ts G,
\ee
of $G$ on the cotangent bundle, induced by the right action of $G$
on itself.
The associated momentum map
\be
\zm :\ts G\simeq\G^\s\ti G\lra\G^\s
\ee
It is a Poisson map with respect to the natural Poisson
structure on $\G^\s$ (see e.g. \cite{AG1, LM}).
\par
Now, we consider any differential 1-form $\zh$ on $\G^\s$ which is
annihilated by the natural Poisson structure $\zL_{\G^*}$
on $\G^\s$ associated with the Lie bracket. Such form will be
called a {\em Casimir form}. We define the vector field $\zG_\zh=
\zL_0(\zm^\s(\zh))$. Then, the corresponding dynamical system can
be written as (for the proof we refer to the general case
described in Theorem 1)
\bea
g^{-1}\dot g&=&\zh(g,p)=\zh(p),\\
\dot p&=&0,
\eea
since $\zw_0=\D (<p,g^{-1}\D g>)$ (cf. \cite{AG1}). Here we
interpret the covector $\zh(p)$ on $\G^\s$ as a vector of $\G$.
Again, our system can be integrated by quadratures, because on each level set,
obtained by fixing $p$'s in $\G^\s$, our dynamical system coincides with a
one-parameter group of the action of $G$ on that particular level set.
\par
We give a familiar example: the rigid rotator and its generalizations.
\begin{ex}
In the case of $G={\rm SO}(3)$ the (right) momentum map
\be
\zm:\ts{\rm SO(3)}\lra{\rm so}(3)^\s
\ee
is a Poisson map onto ${\rm so}(3)^\s$ with the linear Poisson structure
\be
\zL_{{\rm so}(3)^*}=\ze_{ijk}p_i\pa_{p_j}\ot\pa_{p_k}.
\ee
Casimir 1-forms for $\zL_{so(3)^*}$ read $\zh=F\D H_0,$ where
$H_0=\sum p_i^2/2$ is the `free Hamiltonian' and $F=F(p)$ is an arbitrary
function. Clearly, $F\D H_0$ is not a closed form in general, but $(p_i)$ are
first integrals for the dynamical system $\zG_\zh=\zL_0(\zm^\s(\zh))$.
It is easy to see that
\be
\zG_\zh=F(p)\zG_0=F(p)p_i\widehat{X_i},
\ee
where $\widehat{X_i}$ are left-invariant vector fields on ${\rm SO}(3)$,
corresponding to the basis $(X_i)$ of ${\rm so}(3)$ identified with $(\D
p_i)$.
Here we used the identification $\ts{\rm SO}(3)\simeq {\rm SO}(3)\ti{\rm
so}(3)^\s$ given by the momentum map $\zm$.
In other words, the dynamics is given by
\bea
\dot p_i&=&0\\
g^{-1}\dot g&=&F(p)p_iX_i\in{\rm so}(3)
\eea
and it is completely integrable, since it reduces to left-invariant dynamics
on ${\rm SO}(3)$ for every value of $p$.
We recognize the usual isotropic rigid rotator, when $F(p)=1.$
\end{ex}
We can generalize our construction once more, replacing the cotangent bundle
$\ts G$ by its deformation, namely a group double $D(G,\zL_G)$
associated with a Lie-Poisson structure $\zL_G$ on $G$ (see e.g.
\cite{AG2, Lu1}). This double, denoted simply by $D$,
carry on a natural Poisson tensor $\zL^+_D$ which is non-degenerate on the
open-dense subset $D^+=G\cdot G^\s\cap G^\s\cdot G$ of $D$
(here $G^\s\subset D$ is the dual group of $G$ with respect to $\zL_G$).
We refer to $D$ as being {\it complete} if $D^+=D$.
Identifying $D$ with $G\ti G^\s$ if $D$ is complete (or $D^+$ with an open
submanifold of $G\ti G^\s$ in general case; we assume completeness for
simplicity) via the group product,
we can write $\zL^+_D$ in `coordinates' $(g,u)\in G\ti G^\s$ in the form
\be\label{structure}
\zL^+_D(g,u)=\zL_G(g)+\zL_{G^*}(u)-X_i^l(g)\we Y^r_i(u),
\ee
where $X^l_i$ and $Y^r_i$ are, respectively, the left- and right-invariant
vector fields on $G$ and $G^\s$ relative to dual bases $X_i$ and $Y_i$ in
the Lie algebras $\G$ and $\G^\s$, and where $\zL_G$ and $\zL_{G^*}$ are
the corresponding Lie-Poisson tensors on $G$ and $G^\s$ (see \cite{Lu1, AG2}).
It is clear now that the projections $\zm_{G^*}$ and $\zm_G$
of $(D,\zL^+_D)$ onto $(G,\zL_G)$ and $(G^\s,\zL_{G^*})$, respectively,
are Poisson maps. Note that we get the cotangent bundle
$(D,\zL^+_D)=(\ts G,\zL_0)$ if we put $\zL_G=0$.
\par
The group $G$ acts on $(D,\zL_D^+)$ by left translations which, in general
are not canonical transformations. This is, however, a Poisson action with
respect to the inner Poisson structure $\zL_G$ on $G$, which is sufficient
to develop the momentum map reduction theory (see \cite{Lu2}).
For our purposes, let us take a Casimir 1-form $\zh$ for $\zL_{G^*}$, i.e.
$\zL_{G^*}(\zh)=0.$
By means of the momentum map $\zm_{G^*}:D\lra G^\s$, we define the
vector field on $D$:
\be
\zG_\zh=\zL_D^+(\zm_{G^*}^\s(\zh)).
\ee
In `coordinates' $(g,u)$, due to the fact that $\zh$ is a Casimir, we get
\be
\zG_\zh(g,u)=<Y^r_i,\zh>(u)X^l_i(g),
\ee
so that $\zG_\zh$ is associated with the `Legendre map'
\be
L_\zh:D\simeq G\ti G^\s\lra TG\simeq G\ti\G,\quad L_\zh(g,u)=
<Y^r_i,\zh>(u)X_i,
\ee
which can be viewed also as a map $L_\zh :G^\s\lra\G.$
Thus we get the following.
\begin{Theo}
The dynamics $\zG_\zh$ on the group double $D(G,\zL_G)$, associated with a
1-form $\zh$ which is a Casimir for the Lie-Poisson structure $\zL_{G^*}$
on the dual group, is given by the system of equations
\bea
\label{3} \dot u&=&0,\\
\label{4} g^{-1}\dot g&=&<Y^r_i,\zh>(u)X_i\in\G,
\eea
and is therefore completely integrable by quadratures.
\end{Theo}
\begin{ex}
We consider now the double Lie group $D=SL(2,\Bbb C)$ with
$G=SU(2)$ and $G^\s=SB(2,\Bbb C)$ (see e.g. \cite{AG2}.
We will write the elements as  follows:
$$
D=SL(2, \C)\ni a = \left(\begin{array}{cc}
z_1 & z_2 \\ z_3 & z_4\end{array}
\right) ,
     \quad {\rm where }\quad z_i\in \C, \quad z_1z_4-z_2z_3 = 1,
$$
$$
G=SU(2)\ni g = \left(\begin{array}{cc}
 \za & -\bar\zn \\ \zn & \bar \za\end{array}
 \right) ,
     \quad {\rm where }\quad\za, \zn \in  \C, \quad |\za|^2+|\zn|^2=1,
$$
$$
G^\s=SB(2, \C)\ni u =\left(\begin{array}{cc} r& \zg \\ 0 & r^{-1}\end{array}
\right) ,
     \quad {\rm where }\quad r>0,  \zg\in  \C.
$$
The Poisson structure $\zL_{SL(2, \C)}$ is the restriction of
the following quadratic (real) Poisson brackets on $\C^4$:
\be\label{bra}\begin{array}{rclcrcl}
\{z_1,z_2\} &=& -\frac{i}{2} z_1z_2 && \{z_2,z_3\}& =& iz_1z_4 \\
\{z_1,z_3\} &=& \frac{i}{2} z_1z_3 && \{z_2,z_4\}& =&\frac{i}{2}z_2z_4\\
\{z_1,z_4\} &=& 0 &&               \{z_3,z_4\} &=&-\frac{i}{2}z_3z_4\\
\{z_1,\bar z_1\}& =& -\frac{i}{2} |z_1|^2 -i|z_3|^2 &&
     \{z_2,\bar z_2\}& =&-\frac{i}{2} |z_2|^2 -i|z_1|^2-i|z_4|^2  \\
\{z_3,\bar z_3\}& =& -\frac{i}{2} |z_3|^2 &&
     \{z_4,\bar z_4\}& = &-\frac{i}{2}|z_4|^2 -i|z_3|^2\\
\{z_1,\bar z_2\}& =& -iz_3\bar z_4 &&
     \{z_2,\bar z_3\}& =&\frac{i}{2} z_2\bar z_3 \\
\{z_1,\bar z_3\} &=& 0 && \{z_2,\bar z_4\} &=& -iz_1\bar z_3\\
\{z_1,\bar z_4\} &=& \frac{i}{2} z_1\bar z_4 &&\{z_3,\bar z_4\}& =& 0 .
\end{array}
\ee
The lacking commutators may be obtained from this list if we remember
that the Poisson bracket is real, e.g.,
$\{\bar z_i,\bar z_j\} = \overline{\{z_i,z_j\}}$.
One can then check that, indeed, $\det$ and $\overline{\det}$ are
Casimir functions, and that $z_1\leftrightarrow z_4$, $z_2\mapsto -z_2$, and
$z_3\mapsto -z_3$
defines a symmetry of the bracket associated to the inverse
$a\mapsto a^{-1}$ in $SL(2,\C)$.

Our double group is complete, since we have the following unique
(Iwasawa) decompositions:
$$
 SL(2,\C) \simeq SU(2).SB(2, \C),
   \  {\rm where }\quad s=\frac1{\sqrt{|z_1|^2+|z_3|^2}},
$$
\be
\left(\begin{array}{cc} z_1 & z_2\\ z_3 & z_4\end{array}\right)
     = \left(\begin{array}{cc} sz_1 & -s\bar z_3 \\
          sz_3 & s\bar z_1 \end{array}\right)\left(
     \begin{array}{cc}1/s & s(\bar z_1z_2+\bar z_3z_4) \\ 0 &
s \end{array}\right)
\ee
$$
 SL(2,\C) \simeq SB(2,\C).SU(2),
   \  {\rm where }\quad t=\frac1{\sqrt{|z_3|^2+|z_4|^2}},
$$
\be
\left(\begin{array}{cc} z_1 & z_2\\ z_3 & z_4\end{array}\right)
=\left( \begin{array}{cc} t & t(z_1\bar z_3+z_2\bar z_4) \\
          0 & 1/t \end{array}\right)\left(
     \begin{array}{cc} t\bar z_4 & -t\bar z_3 \\ tz_3 &
tz_4 \end{array}\right).
\ee
Hence, the bracket $\{\quad,\quad\}$ is globally symplectic on
$SL(2, \C)$. This bracket is projectable on the subgroups
$SU(2)$ and $SB(2, \C)$, and for the `left trivialization'
$SL(2, \C)=SU(2).SB(2, \C)$ it gives us the Poisson Lie
brackets on $SU(2)$:
\be\begin{array}{rclcrcl}
\{\za,\bar\za\} &=& -i|\zn|^2 && \{\nu,\bar\nu\} &=& 0 \\
\{\za,\nu\}& =& \frac{i}{2} \za\nu &&
     \{\bar\za,\bar\nu\}& =& -\frac{i}{2}\bar\za\bar\nu\\
\{\za,\bar\nu\} &=& \frac{i}{2} \za\bar\nu  &&
     \{\bar\za,\nu\}& =&-\frac{i}{2}\bar\za\nu,
\end{array}
\ee
and on $SB(2, \C)$:
\be
\{\zg,r\} = \frac{i}{2} \zg r,\qquad
     \{\bar\zg,\zg\} = i(r^2-r^{-2}).
\ee
The `interaction' between $SU(2)$ and $SB(2, \C)$ is described by
\be\label{bra1}\begin{array}{rclcrcl}
\{\zn,\zg\}&=&-\frac{i}{4}\zn\zg-i\bar\za r^{-1}&&
\{\za,\zg\}&=&-\frac{i}{4}\za\zg+i\bar\zn r^{-1}\\
\{\bar\zn,\zg\}&=&\frac{i}{4}\bar\zn\zg&&\{\bar\za,\zg\}&=&\frac{i}{4}
\bar\za\zg\\
\{\zn,r\}&=&-\frac{i}{4}\zn r&&\{\za,r\}&=&-\frac{i}{4}\za r,
\end{array}
\ee
where the lacking commutations can be derived, due to the fact that the
bracket
is real.
\par
One can check that the Casimir 1-forms for $\zL_{SB(2, \C)}$
read $\zh=F\D H_0$, where $F=F(\zg,r)$ is an arbitrary function on $SB(2,\C)$
and
$$H_0=\frac{1}{2}Tr(aa^*)=\frac{1}{2}\sum|z_i|^2=\frac{1}{2}(|\zg|^2+
r^2+r^{-2})
$$
is the `free' Hamiltonian.
For the dynamics $\zG_\zh$ on $\SL$ induced by $\zh$,
we calculate (using \ref{bra}) that
\be\label{5}
\begin{array}{lcrclcr}
\dot z_1&=&-\frac{i}{2}F(H_0z_1-\bar z_4)&&
\dot z_2&=&-\frac{i}{2}F(H_0z_2+\bar z_3)\\
\dot z_3&=&-\frac{i}{2}F(H_0z_3+\bar z_2)&&
\dot z_4&=&-\frac{i}{2}F(H_0z_4-\bar z_1).
\end{array}
\ee
Since $F$ and $H_0$ are constants of the motion, it is clear that the system
is completely integrable. For $F=1$ this system was considered in \cite{MSS,
Za}.
In variables $g\in\SU$ and $u\in\SB$, we have
\be
\dot u=\left(\begin{array}{cc}
\dot r&\dot\zg\\0&-\dot rr^{-2}
\end{array}\right)
\ee
and
\be
g^{-1}\dot g=-\frac{i}{2}F(H_0I+J\bar uJu^{-1})=const\in\su,
\ee
where $J=\left(\begin{array}{cc}0&-1\\1&0\end{array}\right)$.
In a more transparent form
\be\label{per}
g^{-1}\dot g=-\frac{i}{4}F\left(\begin{array}{cc}
r^2-r^{-2}+|\zg|^2 & 2\zg r^{-1}\\
2\bar\zg r^{-1} & r^{-2}-r^2-|\zg|^2
\end{array}\right)\in\su.
\ee
It follows that we get a `free motion' on $\SU$ along trajectories of
left-invariant vector fields corresponding to
\be
L_\zh=-\frac{i}{2}F(u)(H_0(u)I+J\bar uJu^{-1})\in\su
\ee
(cf. \cite{MSS, Za}). The mapping
\be
\SB\ni u\mapsto L_\zh(u)=-\frac{i}{2}F(u)(H_0(u)I+J\bar uJu^{-1})\in\su
\ee
is a sort of a `Legendre transform', transforming momenta from $\SB$ into
velocities from $\su$.
It is easy to see that, e.g. in the case $F=1$, $L_\zh$ is
invertible and the momenta corresponding to the velocity
\be-\frac{i}{2}\left(\begin{array}{cc}
s & w \\ \bar w & -s \end{array}\right)\in\su,\quad s\in\R, w\in\C,
\ee
are
\be r=s+\sqrt{s^2+|w|^2+1},\quad \zg=rw.
\ee
\end{ex}
\begin{ex}
Since the dual subgroups in the group double play entirely symmetric role, let
us consider now $\SU$ to be the set of momenta and $\SB$ to be the
configuration space.
\par
The Lie-Poisson structure $\zL_{\SU}$ admits, however no global Casimir
function, so that we cannot produce a globally Hamiltonian system on $\SL$ by
means of the momentum map $\zm_{\SU}:\SL\lra\SU,$ though we easily find out
that
\be
\zh=iF(\za,\zn)(\zn\D\bar\zn-\bar\zn\D\zn)
\ee
is a general Casimir 1-form (which is real for real functions $F$).
The equations of the dynamical system $\zG_\zh$ read
\bea
\dot \za &=& 0\\
\dot \zn &=& 0\\
\dot r &=& -\frac{1}{2}F(\za,\zn)|\zn|^2\zg\\
\dot \zg &=& -F(\za,\zn)\left(\frac{1}{2}|\zn|^2\zg+\frac{i}{r}\za(\Re(\zn)-
\Im(\zn))\right).
\eea
In other words, $\za$ and $\zn$ are constant of the motion and
\be
\dot uu^{-1}=L_\zh(\za,\zn)=-\frac{1}{2}F(\za,\zn)\left(
\begin{array}{cc}
|\zn|^2 & 2i\za(\Re(\zn)-\Im(\zn))\\ 0 & -|\zn|^2
\end{array}\right)\in\Sb
\ee
is time independent.
The `Legendre map' $L_\zh:\SU\lra\Sb$ is never bijective, since our set of
momenta $\SU$ is a compact manifold, so that `admissible' velocities form a
compact subset of $\Sb$.
\end{ex}

\setcounter{equation}{0}
\section{A further generalization}
We have seen that if we concentrate on the possibility of integrating
our system by quadratures, then we can do without the requirement
that the system is Hamiltonian.
\par
By considering again the equations of motion in action-angle variables, we
have, classically,
\bea
\dot I^k&=&0,\\
\dot \zf_k&=&\zn^k(I).
\eea
Clearly, if we have
\bea
\nn\dot I^k&=&F_k(I),\\
\label{6}\dot \zf&=&A_k^j(I)\zf_j,
\eea
and we are able to integrate the first equation by quadratures, we again have
the possibility to integrate by quadratures the system (\ref{6}), if only
the matrices $(A^j_k(I(t)))$ commute:
\be
\zf(t)={\rm exp}\left(\int_0^tA(I(s))\D s\right)\zf_0.
\ee
Of course, because $\zf_k$ are discontinues functions on the torus, we have
to
be more careful here. We show, however, how this idea works for double groups.
In the case when the 1-form $\zh$ on $G^\s$ is not a Casimir 1-form for the
Lie-Poisson structure $\zL_{G^*}$, we get, in view of (\ref{structure}),
\be
\zG_\zh(g,u)=<Y^r_i,\zh>(u)X^l_i(g)+\zL_{G^*}(\zh)(u).
\ee
Now, the momenta evolve according to the dynamics $\zL_{G^*}(\zh)$ on $G^\s$
(which can be interpreted, as we will see later, as being associated with an
interaction of the system with an external field) and `control' the evolution
of the field of velocities on $G$ (being left-invariant for a fixed time) by
a `variation of constants'.
Let us summarize our observations in the following.
\begin{Theo}
The vector field  $\zG_\zh$ on the double group $D(G,\zL_G)$,
associated with a 1-form $\zh$ on $G^\s$,
defines the following dynamics
\bea
\label{7} \dot u&=&\zL_{G^*}(\zh)(u),\\
\label{8} g^{-1}\dot g&=&<Y^r_i,\za>(u)X_i\in\G,
\eea
and is therefore completely integrable, if only we are able to integrate the
equation \ref{7} and $<Y^r_i,\zh>(u(t))X_i$ lie in a commutative subalgebra of
$\G$ for all $t$.
\end{Theo}

\begin{ex}
Let us take now $\SU$ for momenta and consider the Hamiltonian
$H=\frac{1}{2}|\zn|^2.$ The 1-form $\zh=\D H$ is exact, but not a Casimir
1-form for $\zL_{\SU}$. The dynamical system $\zG_\zh$ on $\SL$ induces
the following dynamics of momenta
\bea
\dot \zn &=& 0,\\
\dot \za &=& \frac{i}{2}|\zn|^2\za .
\eea
This is the dynamics described in \cite{LMS}.
Additionally, we get from \ref{bra1}:
\bea
\dot r&=&0,\\
\dot\zg &=& -\frac{i}{2}\bar\za(\Re(\zn)+\Im(\zn))r^{-1}.
\eea
Hence,
$\zn(t)=\zn_0$, $\za(t)=\za_0\exp(\frac{i}{2}|\zn|^2t)$, and
\be
\dot uu^{-1}= -\frac{i}{2}\bar\za(\Re(\zn)+\Im(\zn))\left(\begin{array}{cc}
0 & 1   \\
0 & 0
\end{array}\right)\in\Sb.
\ee
Here, the velocity of a particle is rotating around 0 with the radius
proportional to the momentum $\Re(\zn_0)+\Im(\zn_0)$ and the frequency
proportional to the energy $H=\frac{1}{2}|\zn_0|^2$.
The velocities stay, however, in a commutative subalgebra of the
unipotent Lie algebra $\Sb$, so that
\be
u(t)=\left(\begin{array}{cc}
1 & \za_0\frac{\Re(\zn_0)+\Im(\zn_0)}{|\zn_0|^2}
(\exp(-\frac{i}{2}|\zn_0|^2t)-1)\\
0 & 1
\end{array}\right)u_0.
\ee
\end{ex}
Let us end up with an example which shows that we can actually weaken the
assumptions of Theorem 2. In fact, it is sufficient to assume that
\be
g^{-1}\dot g(t)=\exp(tX)A(t)\exp(-tX)
\ee
for some $A(t),\, X\in\G,$ such that $X+A(t)$ lie in a commutative subalgebra
of $\G$ for all $t$ (e.g. $A(t)=const$), to assure that (\ref{8}) is
integrable
by quadratures. Indeed, in the new variable
\be g_1(t)=\exp(-tX)g(t)\exp(tX)
\ee
the equation (\ref{8}) reads
\be
\dot g_1(t)=g_1(t)(X+A(t))-Xg_1(t)
\ee
and, since the right- and the left-multiplications commute, we easily find
that
\be
g(t)=g_0\exp\left(tX+\int^t_0A(s)\D s\right)\exp(-tX).
\ee
This procedure is similar to what is known as the Dirac interaction picture in
the quantum evolution.
\begin{ex}
For our group double $\SL$, let us take the 1-form $\zh=F(r)\D H_0+\zl\D r,$
$r\in\R$ on $\SB$ which is a perturbation of \ref{per}.
For the dynamics of momenta, we get
\bea
\dot r&=&0\\
\dot \zg&=&-\frac{i}{2}\zl r\zg
\eea
which can be easily integrated:
\bea
r(t)&=&r_0\\
\zg(t)&=&\zg_0\exp(-\frac{i}{2}\zl rt).
\eea
In particular, $r=const$ and $|\zg|=const$. We have also
\be\label{eq}
g^{-1}\dot g=-\frac{i}{4}\left(\begin{array}{cc}
F(r)(r^2-r^{-2}+|\zg|^2)-\zl r & 2F(r)\zg r^{-1}\\
2F(r)\bar\zg r^{-1} & F(r)(r^{-2}-r^2-|\zg|^2)+\zl r
\end{array}\right) .
\ee
The velocities are no longer constant and rotate around the vector
\be\label{A}
C_0=(F(r_0)(r_0^2-r_0^{-2}+|\zg_0|^2)-\zl r_0)\left(
\begin{array}{cc}-\frac{i}{4} & 0 \\ 0 & \frac{i}{4}
\end{array}\right)\in\su.
\ee
It can be interpreted as an effect of an interaction of moving charged
particle with an external magnetic field, corresponding to the perturbation of
the `free system'.
Since, as easily seen,
\be
g^{-1}\dot g(t)=\exp(tX)A_0\exp(-tX),
\ee
with
\be
X=\left(\begin{array}{cc}
-\frac{i}{4}\zl r_0 & 0\\
0 & \frac{i}{4}\zl r_0
\end{array}\right)
\ee
and
$$
A_0=-\frac{i}{4}\left(\begin{array}{cc}
F(r_0)(r_0^2-r_0^{-2}+|\zg_0|^2)-\zl r_0 & 2F(r_0)\zg_0r_0^{-1}\\
2F(r_0)\bar\zg_0r_0^{-1} & F(r_0)(r_0^{-2}-r_0^{2}-|\zg_0|^2)+\zl r_0
\end{array}\right),
$$
we can easily integrate (\ref{eq}):
\bea\nn
g(t)&=&g_0\exp\left(\frac{-iF(r_0)}{4}t\left(\begin{array}{cc}
r_0^2-r_0^{-2}+|\zg_0|^2 & 2\zg_0 r_0^{-1}\\
2\bar\zg_0 r_0^{-1} & r_0^{-2}-r_0^2-|\zg_0|^2
\end{array}\right)\right)\ti\\
&&\left(\begin{array}{cc}
\exp\left(-\frac{i}{4}\zl r_0t\right) & 0\\
0 & \exp\left(\frac{i}{4}\zl r_0t\right)
\end{array}\right).
\eea
\end{ex}
\section{Final comments}
Our identification of systems which can be integrated by quadratures with
one-parameter subgroups of some Lie group $G$ acting on the carrier space of
the system gives us the possibility to dispose of the requirement that the
system is Hamiltonian. We have to notice however that in some cases our system
will turn out to be Hamiltonian with respect to a different symplectic
structure (still invariant under the action of $G$) on the phase space.
Indeed, if we consider $\zh=\zn_i\D I^i$ and make the assumption that
\be
\D\zn_1\we\dots\we\D\zn_n\ne 0
\ee
on some open-dense submanifold $N$, we can define on $N$ the symplectic
structure $\zw_N=\sum_i\D\zn_i\we\D\zvf_i$ and $\zG=\zn_i\pa_{\zvf_i}$ will be
associated with the Hamiltonian $H=\frac{1}{2}\sum_i(\zn_i)^2.$
In the general situation this procedure does not apply any more.
\par
In any case, the importance of the Liouville-Arnold theorem relies on the fact
that, in the hypothesis of the theorem, we can find the group $G$  and
its action on the manifold, and then show that our starting system is
conjugated to the one written in the introduction (\ref{x},\ref{y})
in terms of the action-angle variables.
Our generalization is much more in the spirit of the Lie-Sheffers theorem
\cite{LS} and it
consists of splitting our system along the orbits of an action of a Lie group
and a `transverse' component (which is either zero, or linear), so that the
integration can be achieved easily.
\par
What seems to us relevant is that completely integrable systems (in the
Liouville sense) are a part of this more general scheme. In particular, we
have
shown that we can use the action of the group by non-canonical
transformations,
so that systems on group doubles can be cast in this generalization.
\par
We are confident that this approach may be useful to quantize group doubles in
the geometric quantization setting. Many of these questions are currently
being
investigated.

\end{document}